\documentclass[a4paper,12pt]{amsart}

\usepackage[english]{babel}
\usepackage{amsmath, enumerate, amsfonts, amssymb}%, amsthm}
\theoremstyle{plain}
\newtheorem{theo}{Theorem}[section]
\newtheorem{theoI}{Theorem}
\newtheorem{prop}[theo]{Proposition}

\newtheorem{lem}[theo]{Lemma}
\newtheorem{cor}[theo]{Corollary}
\newtheorem{corI}[theoI]{Corollary}
\newtheorem{clai}[theo]{Claim}

\theoremstyle{definition}
\newtheorem{rem}[theo]{Remark}
\newtheorem*{rem*}{Remark}

\newcommand{\R}{\mathbb{R}}
\newcommand{\C}{\mathbb{C}}
\newcommand{\N}{\mathbb{N}}

\newcommand{\Q}{\mathbb{Q}}
\renewcommand{\k}{\mathbf{k}}

\newcommand{\CC}{\mathcal{C}}

\newcommand{\Frac}{\mathrm{Frac}}

\newcommand{\QQ}{\mathcal{Q}}
\newcommand{\modQ}{{\mathrm{mod}\QQ}}
\newcommand{\PP}{\mathcal{P}}

\newcommand{\spec}{\mathrm{Spec}}

\newcommand{\sdt}{\langle s, \partial_t \rangle}

\renewcommand{\O}{\mathcal{O}} % faisceau des fonctions holomorphes

\newcommand{\D}{\mathcal{D}} % analytic notations

\newcommand{\FD}{\hat{\mathcal{D}}} % formal notations

\newcommand{\dxi}{\partial _{x_i}}
\newcommand{\dtj}{\partial _{t_j}}
\newcommand{\dx}[1]{\partial _{x_{#1}}}
\newcommand{\dt}[1]{\partial_{t_{#1}}}
\newcommand{\ddx}{\partial _x}
\newcommand{\ddt}{\partial_t}

\newcommand{\dxsur}[2]{\frac{\partial {#1}}{\partial x_{#2}}}

\newcommand{\B}{\mathcal{B}}
\newcommand{\ord}{\mathrm{ord}} % order
\newcommand{\gr}{\mathrm{gr}} % graded ring
\newcommand{\lm}{\mathrm{lm}} % leading monomial \lc*\lt
\newcommand{\E}{\mathcal{E}}
\title[Results on Bernstein-Sato polynomials]
{Some results on Bernstein-Sato polynomials
for parametric analytic functions}
\author{Rouchdi Bahloul}
\address{Laboratoire de Math\'ematiques,
Universit\'e de Versailles St-Quentin-en-Yvelines,
45 avenue des Etats-Unis - B\^atiment Fermat,
78035 Versailles,
France}

\begin{document}

%\maketitle

\begin{abstract}
This is the second part of a work dedicated to the study of
Bernstein-Sato polynomials for several analytic functions depending on
parameters. In this part, we give constructive results generalizing
previous ones obtained by the author in the case of one function.
We also make an extensive study of an example for which we
give an expression of a generic (and under some conditions,
a relative) Bernstein-Sato polynomial.
\end{abstract}


\maketitle

%\tableofcontents

Let $X \subset \C^n$ and $Y\subset \C^m$ be compact polydiscs centered
at the origin, $Z=X \times Y$ and $f=(f_1,\ldots,f_p)$ ($p\ge 2$) an
analytic map from $X$ to $\C^p$. We are interested in the study of
Bernstein-Sato polynomials of $f(x,y_0)$ when $y_0$ moves through $Y$.
Our work is related to the notion of generic
Bernstein-Sato polynomials as in Brian\c{c}on et al. \cite{bgm} (for
$p=1$) and Biosca \cite{biosca1}. Herein we shall adopt a more
constructive method as in Bahloul \cite{jmsj} (where the case $p=1$
was treated), based on the first part \cite{part1} and Bahloul
\cite{compos}.

Our goal is to give analogous results to \cite{jmsj}. However, since
the construction in \cite{compos} is entirely algorithmic only when
$p=2$, a part of the results herein shall be shown only for $p=2$. It
would be a nice result if one could wholly achieve \cite{compos} in an
algorithmic way (here ``algorithmic'' means ``in an infinite way'').
Note that a similar question was treated in the case of polynomials
$f_j$ in Bahloul \cite{procJap} with direct methods while
constructive methods were used in Leykin \cite{leykin} (for $p=1$) and
Brian\c{c}on, Maisonobe \cite{bm02} (for $p\ge 1$).

{\bf Note.}
If $\O_{\C^{n+m}}$ denotes the sheaf of analytic functions on
$\C^{n+m}$, we shall identify $\O_Z$ with the germ
$\O_{\C^{n+m},0}$. Sometimes, we will reduce $Z$ without an explicit
mention so that $\O_Z$ shall be identified with the set
$\O_{\C^{n+m}}(U)$ of sections of $\O_{\C^{n+m}}$ on an open (poly)disc
$0 \in U \subset Z$.

\section{Main results}

$\D_{Z/Y}$ denotes the ring of relative differential operators. It is
the subring of $\D_Z$ made of elements without derivations
$\partial_{y_i}$. Let us write $s=(s_1,\ldots,s_p)$
and $\ddt=(\dt{1}, \ldots, \dt{p})$. Following \cite{bm02}, define
$\C\sdt$ as the algebra $\C[s, \ddt]$ with the relations $\dtj s_j=
s_j \dtj -\dtj$ ($j=1, \ldots, p$) and set $\D_Z \sdt=\D_Z \otimes
\C\sdt$. If $t=(t_1, \ldots,t_p)$ are new indeterminates, the
identification $s_j=-\dtj t_j$ gives the inclusions of rings: $\D_Z[s]
\subset \D_Z\sdt \subset \D_{Z \times \C^p}$. This identification
comes from the fact that the free $\O_Z[1/F, s]$-module $\O_Z[1/F, s]
\cdot f^s$ (here $F=f_1 \ldots f_p$ and $f^s=f_1^{s_1} \cdots
f_p^{s_p}$) is a $\D_{Z \times \C^p}$-module and the action of $s_j$
coincides with that of $-\dt{j} t_j$ (see Malgrange
\cite{malgrange}).

For a given set of (germs at $0$ of) analytic functions $g=(g_1,
\ldots, g_p)$ on $X$, $\B(g)$ shall denote the ideal of Bernstein-Sato
of $g$ (at $x=0$):
it is the set of $b(s)\in \C[s]$ satisfying $b(s) g^s \in
\D_X[s] g^{s+1}$ (here $g^{s+1}:=\prod g_j^{s_j+1}$). This ideal is
not zero (Sabbah \cite{sabbah}) and in fact it contains a polynomial
of the form $\prod (l_1 s_1 +\cdots + l_p s_p +a)$ with $l_j
\in \N$ and $a \in \Q_{>0}$ (Gyoja \cite{gyoja}).

\begin{rem}\label{facteurs}
If for any $j=1,\ldots ,p$, $g_j^{-1}(0) \nsubseteq
\bigcup_{k \ne j} g_k^{-1}(0)$ then $\B(g) \subset \C[s]
\cdot \prod_j (s_j+1)$.
Indeed, it suffices to specialize $s_j=-1$ in a functional
equation.
\end{rem}

When the $g_j$ are in $\k[[x]]$ for some field $\k$ (of
characteristic $0$) we can also consider $\B(g) \subset \k[s]$ the
ideal defined by the same relation where we replace $\D_X$ by
$\FD_x(\k)=\k[[x]][\ddx]$. It is well known \cite{bri-mais}
that given $g \in (\O_X)^p \subset \C[[x]]$,
the formal Bernstein-Sato ideal coincides with the analytic one.

For $g \in (\k[[x]])^p$ it is still an open question whether or not
$\B(g)$ is zero. We know it is not zero only when $p=1$
(Bj\"ork \cite{bjork}).

Let us come back to our situation.
We retain the notations of part 1 \cite{part1}. Set $\CC=\O_Y$ and
$\QQ \in \spec(\CC)$.
Each $f_j$ is viewed as an element of
$\CC[[x]]$ and we consider the specialization $(f_j)_\QQ \in
\Frac(\CC/\QQ)[[x]]$ to $\QQ$ so that $\B((f)_\QQ)$ is an ideal of
$\Frac(\CC/\QQ)[s]$.

\begin{theoI}\label{theoI1}
For any $b(s) \in \C[s]$, the following conditions are equivalent:
\begin{itemize}
\item[(i)]
$b(s) \in \B((f)_\QQ)$.
\item[(ii)]
$\exists h(x,y) \in \O_Z$ with $h(0,y) \notin \QQ$ such that
%\begin{equation}\label{eq:faible}
\[h(x,y) b(s) f^s \in \D_{Z/Y}[s] f^{s+1} + \QQ \D_{Z/Y} \sdt f^s.\]
%\end{equation}
\item[(iii)]
$\exists c(y) \in \O_Y \smallsetminus \QQ$ such that for any
$y_0 \in V(\QQ) \smallsetminus V(c)$, $b(s) \in \B(f(x,y_0))$.
\end{itemize}
\end{theoI}

\begin{proof}
The proof is similar to that of \cite[Th. 1]{jmsj}.
Let us give its main lines.
For $(ii) \Rightarrow (i)$, it suffices to specialize to $\QQ$,
while $(ii) \Rightarrow (iii)$ is trivial by taking $c(y)=h(0,y)$.

Let us introduce two ideals:
$I_0$ the ideal of $\FD_x(\CC)\sdt=\CC[[x]][\ddx] \otimes \CC\sdt$
generated by the
$s_j+f_j \dtj$, $j=1,\ldots,p$, and the $\dxi+\sum_{j=1}^p
\dxsur{f_j}{i} \dtj$, $i=1,\ldots,n$;
and $I_0'\subset \FD_x(\CC)[s]$ defined as
\[I_0'=(I_0 +\FD_x(\QQ)\sdt) \cap \D_X(\CC)[s]+ \D_X(\CC)[s] \cdot
F.\]

For any $y_0 \in Y$, ${I_0}_{|y=y_0}$
is the annihilator in $\FD_x(\C)\sdt$ of $\prod_j f_j(x,y_0)^{s_j}$
(see e.g. \cite[sect. 4]{jmsj} for $p=1$, the proof for $p\ge 2$ is
the same).

Moreover, we have the following (by using an arbitrary generic
standard basis of $I_0'$): $(I_0')_\QQ$ equals
\[(I_0)_\QQ \cap \FD_x(\Frac(\CC/\QQ))[s]+
\FD_x(\Frac(\CC/\QQ))[s] \cdot (F)_\QQ.\]

As a consequence, for $b \in \C[s]$, $b \in \B((f)_\QQ)$ if and only
if $b\in (I_0')_\QQ$.

Now assume we have (iii).
Consider the division modulo $\QQ$ of $b$ by (a generic standard basis
of) $I_0'$ (see \cite[Prop. 2.2]{part1} and
\cite[Prop. 3.5]{jmsj})
and denote by $R$ the remainder
$\modQ$. It follows that for a generic $y_0 \in V(\QQ)$, $R_{|y=y_0}$
is in ${I_0'}_{|y=y_0}$. This is possible only if $R$ is zero modulo
$\QQ$, thus
\[b \in \FD_x(\CC[c^{-1}])[s] \cdot I_0' + \FD_x(\CC[c^{-1}])[s] \cdot
F\]
for some $c\in \CC \smallsetminus \QQ$.
Specializing this relation to $\QQ$, we get $b=(b)_\QQ \in
(I_0')_\QQ$. Thus (i) is satisfied.

Now assume we have (i), which means that $b \in (I_0')_\QQ$.
Let us consider the division modulo $\QQ$ of $b$ by (an arbitrary
generic standard basis of) $I_0'$.
The remainder is zero modulo $\QQ$, which means that
\[b \in \FD_x(\CC[c^{-1}])[s] \cdot I_0' + \FD_x(\QQ[c^{-1}])[s]\]
for some $c \in \O_Y \smallsetminus \QQ$.
Applying $b$ to $f^s$, we obtain a formal functional equation of
the form
\[b f^s \in \FD_x(\O_Y[c^{-1}])[s] f^{s+1} + \FD_x(\QQ[c^{-1}])
\sdt f^s.\]
We may then pass from the formal to the analytic setting
(following the same arguments as in the last section of \cite{jmsj})
and we get (ii).
\end{proof}

We still don't know whether or not $\B((f)_\QQ) \cap \C[s]$ is zero.

\begin{theoI}\label{theoI2}
Here $p=2$ (see the comments in $\S 2$ of the introduction).
There exists a non zero polynomial $b(s)$ of the form $\prod (l_1 s_1
+\cdots + l_p s_p +a)$ with $l_j \in \N$ and $a \in \Q_{>0}$, that
belongs to $\B((f)_\QQ)$.
\end{theoI}
The proof will be given in section 4.

As a consequence: there exists a finite stratification
$Y=\cup W$ into locally closed subsets $W$ and polynomials $b_W(s)$
of the above form such that for any $y_0 \in W$, $b_W(s) \in
\B(f(x,y_0))$.

Consider the lcm of the $b_W$ and denote it by $b_{\mathrm{comp}}$
then for any $y_0 \in Y$, $b_{\mathrm{comp}}$ is a Bernstein-Sato
polynomial of $f(x,y_0)$.
Here $b_{\mathrm{comp}}$ should be read ``comprehensive Bernstein-Sato
polynomial''. It is clear that any ``relative Bernstein-Sato
polynomial'' is comprehensive but the converse is obviously wrong
since a relative Bernstein-Sato polynomial does not exist in general,
even when $p=1$
(see e.g. \cite{biosca2} for the definition of a relative
Bernstein-Sato polynomial, see also \cite{BM-relative} for general
results on the subject in the hypersurface case).

\begin{corI}\label{corI}
Here $p$ is not necessarily $2$.
Take $n=2$ and suppose that for a generic $y_0$ in $V(\QQ)$,
$f_1(x,y_0), \ldots, f_p(x,y_0)$ are irreducible and pairwise
relatively prime. Take $b(s) \in \C[s]$ then $b(s) \in
\B((f)_\QQ)$ if and only if
there exists $H(y) \in \O_Y \smallsetminus \QQ$ with
\[ H(y) b(s) \in \D_{Z/Y}[s] f^{s+1} + \QQ \D_{Z/Y} \sdt f^s.\]
\end{corI}
This means that $b(s)$ is a ``generic Bernstein-Sato'' polynomial in
the sense of Biosca \cite{biosca1} (notice that in previous works on
this subject, the notion of generic Bernstein-Sato polynomial is
defined only when $\QQ=(0)$, see e.g. loc. cit. and its references).

The assumptions of this corrolary mean that
the relative singular locus $V(\dxsur{F}{1}, \ldots, \dxsur{F}{n}, F)$
projects to $0$ by the projection $X \times Y \to X$ when we restrict
ourself to $X \times U$ and $U$ is a Zariski open set of $V(\QQ)$.

Let us give a:
\begin{proof}[Sketch of Proof of Cor. \ref{corI}]
The ``if'' sense is trivial. Let us prove the converse.
We don't give all the details of the proof for it is analogous
to that of \cite[Cor. 2]{jmsj}.
Denote by $J$ the ideal of $\O_Z$ generated by $F$ and the
$\dxsur{F}{i}$'s. The hypothesis can be rephrased as follows:
\[V(\sqrt{J+\O_Z \cdot \QQ} : h_0) \subset (0) \times V(\QQ)\]
in $Z=X\times Y$, for some $h_0(y) \in \O_Y \smallsetminus \QQ$.
Thus the zero locus of $\sqrt{J+\O_Z \cdot \QQ} : h_0 + \O_Z \cdot h$
is included in the zero set of $h(0,y)$, where $h$ is obtained
from Th. \ref{theoI1}(ii).
As a consequence $H:=(h_0 h(0,y))^k$ is in $J+ 
\O_Z \cdot \QQ +\O_Z \cdot h$, for some $k\in \N$, by using
Hilbert's Nullstellensatz. This $H$ is in $\O_Y \smallsetminus
\QQ$.

Now, for a generic $y_0$ in $V(\QQ)$, $b(s)$ is a Bernstein-Sato
polynomial for $f(x,y_0)$, thus by assumption and Rem. \ref{facteurs},
$\prod_1^p (s_j+1)$ divides $b(s)$. Let us
write $b(s)=\prod_1^p (s_j+1) \cdot \tilde{b}(s)$. For
$i=1, \ldots, n$, we have:
\[
\begin{array}{l}
b(s) \dxsur{F}{i} f^s=\\
\tilde{b}(s) \prod_{k=1}^p(s_k+1) (\sum_{j=1}^p \dxsur{f_j}{i}
     \frac{F}{f_j}) f^s=\\
\tilde{b}(s) \sum_{j=1}^p (\prod_{k\ne j} (s_k+1))
(s_j+1) \dxsur{f_j}{i} \frac{1}{f_j} f^{s+1}=\\
\tilde{b}(s) \sum_{j=1}^p (\prod_{k\ne j} (s_k+1)) \dxi \cdot
f^{s+1}.
%b(s) \dxsur{F}{i} f_1^{s_1} \cdots f_p^{s_p}=\\
%\tilde{b}(s) \prod_{k=1}^p(s_k+1) (\sum_{j=1}^p \dxsur{f_j}{i}
%     \frac{F}{f_j}) f_1^{s_1} \cdots f_p^{s_p}=\\
%\tilde{b}(s) \sum_{j=1}^p (\prod_{k\ne j} (s_k+1))
% (s_j+1) \dxsur{f_j}{i} \frac{1}{f_j} f_1^{s_1+1} \cdots
%f_p^{s_p+1}\\
%=\tilde{b}(s) \sum_{j=1}^p (\prod_{k\ne j} (s_k+1)) \dxi \cdot
%f_1^{s_1+1} \cdots f_p^{s_p+1}
\end{array}
\]
From this equality, and relation (ii) in Th. \ref{theoI1},
we get the desired equation with $H(y)$.
\end{proof}

\section{An example related to \cite{kyushu}}

Let us consider the following example:

$f_1(x,y)=c_1(y) x_1^a+ c_2(y) x_2^b+g_1(x_1,x_2, y)$

$f_2(x,y)=c_3(y) x_1^c+ c_4(y) x_2^d+g_2(x_1,x_2,y)$

\noindent
with $y=(y_1,\ldots,y_m)$, $x=(x_1,x_2)$, $a,b,c,d \in \N_{>0}$.
Here, $c_i \in \O_Y$ and $g_i(x,y) \in \O_Z$.
We assume that $C(y):=\prod_i c_i(y)$ is not zero, and we work
with $\QQ=(0)$ so that $V(\QQ)=Y$.

Consider the weight vectors $\alpha_1=(b,a)$, $\alpha_2=(d,c)$ on
the variables $(x_1,x_2)$, with $bc > ad$. The weight of an element
$g$ in $\O_Z$ for $\alpha_i$, denoted by $\rho_{\alpha_i}(g)$,
is the minimum of the $\alpha_i$-degrees in $x$ of the monomials of
$g$.

We assume that $\rho_{\alpha_i}(g_i) > \rho_{\alpha_i}(f_i)$, $i=1,2$.
As a consequence, for any $y_0$ with $C(y_0)\ne 0$,
it is easy to check that $f_1(x,y_0)$ and $f_2(x,y_0)$ are
irreducible and relatively prime, so $f$ satisfies the assumptions of
Cor. \ref{corI}.

On the other hand, for any $y_0$ with $C(y_0) \ne 0$,
the main result of Bahloul \cite{kyushu} applies.
Put $N_1=2ab+ad-2a-2b$, $N_2=2cd+ad-2c-2d$,
$W_1=\{\deg_{\alpha_1}(z)\}$ (resp. $W_2=\{\deg_{\alpha_2}(z)\}$), $z$
running over the monomials with
$\deg_{\alpha_1}(z) \le N_1 +\rho_{\alpha_1}(f_2)$
(resp. $\deg_{\alpha_2}(z) \le N_2 +\rho_{\alpha_2}(f_1)$), and
$b(s_1,s_2)=(s_1+1)(s_2+1) \prod_{\rho_1 \in W_1} (ab s_1+ad s_2+
a+b+\rho_1) \prod_{\rho_2\in W_2} (ad s_1+cd s_2+c+d+ \rho_2)$.

By \cite[Prop. 1]{kyushu}, for any $y_0\in Y$ with $C(y_0)\ne 0$,
the polynomial $b(s_1,s_2)$ is in $\B((f_1,f_2)(x,y_0))$, that is,
$b$ satisfies Th. \ref{theoI1}(iii). Applying Cor. \ref{corI}, we get:
\[H(y) b(s_1,s_2) f_1^{s_1} f_2^{s_2} \in \D_{Z/Y}[s_1,s_2]
f_1^{s_1+1} f_2^{s_2+1}\]
for some non zero $H(y)\in \O_Y$. This means that $b$ is a generic
Bernstein-Sato polynomial in the usual sense. If we look at the
details of the proof of (iii)$\Rightarrow$(ii) in Th. \ref{theoI1} and
the proof of Cor. \ref{corI}, we notice that the $H(y)$ obtained in
this corollary is of the form $C(y)^k$ for some $k\in \N$, thus:

\begin{prop}
For some $k\in \N$, we have
\[C(y)^k b(s_1,s_2) f_1^{s_1} f_2^{s_2} \in \D_{Z/Y}[s_1,s_2]
f_1^{s_1+1} f_2^{s_2+1}.\]
\end{prop}
As a consequence, if $C(y)$ is invertible (i.e. $C(0) \ne 0$) then $b$
is a relative Bernstein-Sato polynomial. In fact, we have a more
precise statement:

\begin{prop}
Set $C'(y)=c_1(y) c_4(y)$.
If $C'$ is invertible then the polynomial $b$ above is a relative
Bernstein-Sato polynomial:
\[b(s_1,s_2) f_1^{s_1} f_2^{s_2} \in
\D_{Z/Y}[s_1,s_2] f_1^{s_1+1} f_2^{s_2+1}.\]
\end{prop}

It is a direct consequence of the following result.

\begin{clai}
The polynomial $b$ above satisfies:
\[(1+p) {C'}^k b(s_1,s_2) f_1^{s_1} f_2^{s_2} \in
\D_{Z/Y}[s_1,s_2] f_1^{s_1+1} f_2^{s_2+1}\]
for some $p \in \sum_{i=1}^n\O_Z[{C'}^{-1}]\cdot x_i$ and some $k\in
\N$.
\end{clai}

%Notice that from a geometric point of view, the notion of relative
%Bernstein-Sato polynomial is not well known for $p\ge 2$ (see Biosca
%\cite{biosca2}).

\begin{proof}[Proof of the claim]
The proof follows \cite{kyushu}. Let us first review it and then
explain how it can be adapted to our situation.
In \cite{kyushu}, the data are analytic functions $f_1, f_2$
satisfying some conditions. For example, $f_1=x_1^a+x_2^b +g_1(x)$
and $f_2= x_1^c+x_2^d+g_2(x)$ with $\rho_{\alpha_i}(g_i) >
\rho_{\alpha_i}(f_i)$.
We define $\xi_{i_1,i_2}=\prod_{k=0}^{i_1-1}(s_1-k)
\prod_{k=0}^{i_2-1} (s_2-k) f_1^{s_1-i_1} f_2^{s_2-i_2}$ for
$(i_1,i_2) \in \N^2$, $\xi_{-1,0}=f_1^{s_1+1}f_2^{s_2}$,
$\xi_{0,-1}=f_1^{s_1}f_2^{s_2+1}$,
$\xi_{-1,-1}=f_1^{s_1+1}f_2^{s_2+1}$.
Then we attach $\alpha_i$-weights to the elements of $\D_X[s_1,s_2]
\xi_{i_1,i_2}$ (see \cite[Def. 1.3]{kyushu}), by defining
$\rho_{\alpha_i}(\sum_{\beta, k,l} \ddx^\beta s_1^k s_2^l
u_{\beta k l}(x) \xi_{i_1,i_2})$ as the minimum of $\rho_{\alpha_i}
(u_{\beta k l}(x))-i_1 \rho_{\alpha_i}(f_1) -i_2 \rho_{\alpha_i}(f_2)$.

On the other hand, we introduce the ideals of $\O_X$: $I=\langle f_1,
f_2 \rangle$, $I_1=\langle f_1, J \rangle$ and
$I_2=\langle f_2, J \rangle$. Here $J$ is the determinant of the
jacobian matrix of $(f_1,f_2)(x_1,x_2)$.
We show that these ideals have a finite colength lower
than $N_1$ and $N_2$. For this purpose we use divisions and standard
bases settings. The local order used in the divisions is such that
the leading terms of
$f_1$, $f_2$ and $J$ are $x_1^a$, $x_2^d$, $a d x_1^{a-1} x_2^{d-1}$
respectively.

\emph{Step 1.}
The first step of the proof is to show that applying $b$ to
$\xi_{0,0}=f_1^{s_1} f_2^{s_2}$ gives rise to a (finite) sum of
elements $(s_1+1)(s_2+1) P(s) \xi_{i_1, i_2}$ with $\alpha_i$-weight
$>N_i$.

\emph{Step 2.}
By division first by $I$ and then by the $I_i$'s we can go down
from $\D_X[s] \xi_{i_1,i_2}$ to $\D_X[s] \xi_{i_i-1, i_2}$ and
$\D_X[s] \xi_{i_1, i_2-1}$ while the $\alpha_i$-weight is
conserved. This enables an induction on $i_1$ and $i_2$, so that we
can go back to $\xi_{0,0}=f_1^{s_1} f_2^{s_2}$.

Now let us see how the proof of \cite{kyushu} can be adapted to prove
our claim.

We work in a formal setting so that $f_i$ are viewed in $\O_Y[[x]]$.
Step 1 can be done without any problems. In step 2, we shall do
divisions by $I$ (resp. $I_1$, $I_2$). But the
leading terms of $f_1, f_2, J$ are $c_1 x_1^a$, $c_4 x_2^d$,
$c_1 c_4 a d x_1^{a-1} x_2^{d-1}$ (see Propr. 2.5 and the proof of
Aff. 4.1 in \cite{kyushu})
so all the divisions will take place in $\O_Y[{C'}^{-1}][[x]]$.
Therefore, the equation obtained will be of the form:
\[b(s_1,s_2) f_1^{s_1} f_2^{s_2} \in
\FD(\O_Y[{C'}^{-1}])[s_1,s_2] f_1^{s_1+1} f_2^{s_2+1}.\]
As above, we may pass from the formal to the analytic setting
to conclude.
\end{proof}

\section{Recalls and preliminaries on Bernstein-Sato polynomials}

In order to prove Th. \ref{theoI2}, we shall review some results.

For more details, see \cite{compos}.
The system of coordinates $(x,t)$ being fixed,
we denote by $V_j$ the $V$-filtration associated with the hypersurface
$t_j=0$, $j=1,\ldots,p$, on $\D_{X \times \C^p}$ (we can see it as the
natural filtration associated with the weight vector also denoted
$V_j$  where the weight of $t_j$ and $\dt{j}$ are $-1$ and $1$
respectively, and the weight is zero for the other symbols). For
$L=(l_1, \ldots, l_p)$ in $(\R_{\ge 0})^p$,
we denote by $V^L$ the filtration $\sum_j l_j V_j$ and $\gr^L$ the
associated graded ring.

Given $g=(g_1,\ldots,g_p)$ analytic on $X$, we define $I$ the
annihilator of $g^s$ in $\D_{X \times \C^p}$.
% and we put $M=\D_{X \times \C^p}/I$.
The ideal $\B_L$ is then defined as the
set of $c(s) \in \C[s_1, \ldots, s_p]$ with the relation
\[c(s) g^s \in V_{<0}^L(\D_{X \times \C^p}) g^s .\]
Then $b_L=b_{L,g}$ (if it is not zero) is the monic polynomial
$e(\lambda)$ in one variable of the least degree satisfying:
$e(L(s)) \in \B_L$. Here $L(s)=\sum_j l_j s_j$.
This polynomial is not zero (Sabbah \cite{sabbah} and has roots in
$\Q_{<0}$ Gyoja \cite{gyoja}). It can be seen as the Bernstein-Sato
polynomial of $g$ in the ``direction'' $L$. Notice that in the
algebraic case, $b_{(1,\ldots,1)}$ coincides with the $b$-function
considered in Budur et al. \cite{bms}.

Now we can consider the restriction $\E_V(h(I))$ to the space
$\sum_j \R_{\ge 0} V_j$ of the analytic Gr\"obner fan of $I$, for
which we denote by $Sq(\E_V(h(I)))$ the $1$-skeleton. This restriction
leaves in $(\R_{\ge 0})^p$ and this skeleton is in $\N^p$ (because the
Gr\"obner fan is rational).

\begin{theo}[\cite{sabbah} and \cite{compos}]\label{theo:sab-bah}
There exists $\kappa \in \N^p$ such that the polynomial
\[b(s)= \prod_{L \in Sq(\E_V(h(I)))} \prod_{-L(\kappa +
(1,\cdots,1)) <k \le 0 } b_{L,g}(L(s)-k) \]
is Bernstein-Sato polynomial of $g$.
\end{theo}

\begin{rem}\label{rem:kappa}
\begin{itemize}
\item
In Sabbah \cite{sabbah}, the author shows that there exists $\kappa
\in \N^p$ satisfying a certain property, say $(\mathcal{P})$.
Then his shows that if $\kappa$ satisfies $(\mathcal{P})$ then it
satisfies Th. \ref{theo:sab-bah}.
\item
When $p=2$, in \cite{compos}, we construct explicitely some $\kappa$
making $(\mathcal{P})$ true. The construction goes
as follows: Let $L_1, \ldots, L_q \in \N^p$ be such that
$C_{L_i}(h(I))$ are the maximal cones of the (open) fan
$\E_V(h(I))$. Let $G_i$ be the reduced
standard basis of $h(I)$ for an order $\prec_{L_i}^h$
adapted to $L_i$ and define
$\kappa^1$ as the maximum of the $\ord^{V_1}(P)-\ord^{V_1}(
\lm_{\prec_{L_i}^h}(P))$ where $P$ runs over all the element of all
the $G_i$'s. Then $\kappa=(\kappa^1, 0)$ satisfies property
$(\mathcal{P})$. Here $\ord^{V_1}$ means the order with respect to the
filtration $V_1$ (see \cite{compos}).

Notice that this $\kappa$ depends (only) on two monomials of each
element of the standard bases. Thus it depends on a finite number of
monomials.
\end{itemize}
\end{rem}

\begin{lem}\label{lem:B_L}
\begin{itemize}
\item[(a)]
With the identification $s_j=-\dtj t_j$, we have $\B_L(g)= \gr^L(I)
\cap \C[s_1,\ldots,s_p]$.
\item[(b)]
$b_{L,g}(L(s))$ is the monic generator of $\gr^L(I) \cap \C[l_1 s_1+
\cdots + l_p s_p]$.
\end{itemize}
\end{lem}
The lemma is a straightforward consequence of the definitions.

\section{Proof of Theorem \ref{theoI2}}

The proof shall be partially sketched because it uses the same method
as in \cite{jmsj}.

\subsection{Formal algorithm for $b_L$}

Given $L \in (\R_{\ge 0})^p$,
we give an algorithm for computing the polynomial $b_L$
for formal power series $g_1, \ldots, g_p \in \k[[x]]$, $\k$
denotes a field of characteristic $0$.

\begin{lem}
Let $I$ be an ideal in $\FD_{x,t}$. Let $2\le k \le p$ and
$\{j_1, \ldots, j_k\} \subset \{1,\ldots, p\}$ such that
$l_{j} = 0$ iff $j=j_i$ with $1\le i \le k$
(if none of the $l_j$ is $0$ then put $k=1$).
Then $\gr^L(\FD_{x,t})$ is canonically isomorphic to
$\k[[x, t_{j_1}, \ldots, t_{j_{k-1}}]]\langle t_{j_k}, \ldots,
t_{j_p}, \ddt, \ddx \rangle$, where the commutation relations are
obtained from $\FD_{x,t}$ by restriction.
\end{lem}
The proof is straightforward. Thus this graded ring is a subring of
$\FD_{x,t}$ and it can be constructed as in \cite[Section 3]{jmsj}.
Therefore all the results of loc. cit. about (generic) standard
bases apply.

In the following, in order to simplify, we assume $\{j_1, \ldots,
j_{k-1}\}=\{1, \ldots, k-1\}$, i.e. $L=(0, l_k, \ldots, l_p)$.
Consider the following ideals:
\begin{itemize}
\item[(0)]
$I=\sum_{j=1}^p \FD_{x,t} (t_j -g_j) + \sum_{i=1}^n
\FD_{x,t} (\dxi+ \sum_{j=1}^p \dxsur{g_j}{i} \dtj)$.\\
This ideal is the annihilator of $g^s$ in $\FD_{x,t}$.
\item[(1)]
$I_1=\gr^L(I)$ in\\ $\k[[x, t_{j_1}, \ldots, t_{j_{k-1}}]]
\langle t_{j_k}, \ldots, t_{j_p}, \ddt, \ddx \rangle$.
\item[(2)]
$I_2= I_1 \cap \k[[x,t_1,\ldots,t_{k-1}]]
\langle t_k, \ldots, t_p, \dt{k}, \ldots, \dt{p} \rangle$.
It is an elimination of the ``global'' variables $\dx{i}$ and
$\dt{j}$ for $j=1, \ldots, k-1$.
\item[(3)]
We introduce a new indeterminate $\lambda$ and we consider the
ring $\k[[x,t_1,\ldots,t_{k-1}]] \langle
t_k, \ldots, t_p, \dt{k}, \ldots, \dt{p} \rangle [\lambda]$
where the new relations are: $[\lambda, x_i]=[\lambda, t_j]=0$
for any $i$ and any $1 \le j \le k-1$ and for $j \ge k$
$[t_j, \lambda]=l_j t_j$ and $[\dtj, \lambda]=-l_j \dt{j}$.
In other terms, $\lambda$ behaves like $l_k s_k +\cdots +l_p s_p$
where $s_j=- \dtj t_j$.\\
We consider then the previous ideal in this ring and we put\\
$I_3=I_2 \cap \k[[x, t_1, \ldots, t_{k-1}]][\lambda]$.\\
We have eliminated the ``global'' variables $t_j$ and $\dtj$
for $j \ge k$. Notice that now we have a commutative setting.
\item[(4)]
$I_4=I_3 \cap \k[\lambda]$.
We eliminate the ``local'' variables $x_i$ and $t_j$.
\end{itemize}

Lemma \ref{lem:B_L} asserts that the monic
generator of $I_4$ (if it is not zero) is the polynomial
$b_{L,g}$.

The reason why we need to go through step $(3)$ is that we know
how to go from
a given ideal $\mathcal{I} \subset \k[[y_1,\ldots,y_m]][\lambda_1,
\ldots, \lambda_q]$
to the ideal $\mathcal{I} \cap \k[\lambda_1,\ldots,\lambda_q]$ only
when $q=1$. We can
find such an algorithm in \cite[4.1]{jmsj} (which is a variant
of Oaku's \cite[Algo. 4.5]{oaku}).\\

\noindent
\emph{Details for step (1).}
All the steps but step (1) consist in the elimination of global or
local variables.
The elimination of global variables can be done as in
\cite[Prop. 3.8]{jmsj} whereas the local elimination is described
in \cite[4.1]{jmsj}. Let us discuss step (1).
As in \cite{part1}, we consider $h(I)
\subset \FD_{x,t}\langle h \rangle$ (we use $h$ instead of $z$
not to make confusions with \cite{jmsj}) generated by the
degree-homogenization of the elements of $I$. This ideal can
be obtained via a standard basis with respect to an order that
respects the total degree in the $\dxi$'s. Then we compute a
$\prec_L^h$-standard basis $G$ and $G_{|h=1}$ shall be a
system of generators of $\gr^L(I)$.
This is well known, see for example \cite{btaka}.

In conclusion,
we can get $\gr^L(I)$ from $I$ via standard bases computations
for (admissible) orders as in \cite{part1, jmsj}.

\subsection{Here is the proof}

For the proof, we shall use tools from \cite{jmsj}
and the first part \cite{part1}.

Consider $b_{L,((f)_\QQ)}$. We don't know a priori whether or not it
is zero. By \cite[sect. 4-5]{jmsj} applied to the previous
algorithm, we have: for a generic $y_0 \in V(\QQ)$,
$b_{L, f(x,y_0)}$ is constant and
equal to $b_{L,((f)_\QQ)}$. So the latter is not zero and has
rational roots. (Notice that this argument is valid for any
$p\ge 2$).

Now consider the  ideal $I$ in $\D_{Z \times \C^p/Y}$ generated by
the $t_j-f_j$, $j=1,\ldots,p$, and the $\dxi+\sum_{j=1}^p
\dxsur{f_j}{i} \dtj$, $i=1,\ldots,n$.
From \cite{part1}, we know that the analytic Gr\"obner fan of $I$
is constant for a generic $y_0 \in V(\QQ)$, so that the same is true
for its restriction to the space $\sum_{j=1}^p \R_{\ge 0} V_j$, and it
equals $\E_V(h((I)_\QQ))$.
Moreover, it follows from Remark \ref{rem:kappa} (it is here that
we need to assume $p=2$) that the element $\kappa$ obtained from the
Gr\"obner fan (as it is explained in this remark) is generically
constant. This $\kappa$ satisfies property $\mathcal{P}$ for
$f(x,y_0)$ for a generic $y_0 \in V(\QQ)$. This implies that 
the polynomial
\[\prod_{L \in Sq(\E_V(h((I_\QQ))))} \prod_{-L(\kappa +
(1,\cdots,1)) <k \le 0 } b_{L,(f)_\QQ}(L(s)-k) \]
belongs to $\B(f(x,y_0))$ for a generic $y_0 \in V(\QQ)$.
We then apply (iii) in Th. \ref{theoI1} to conclude.

\end{document}